\documentclass{article}
\usepackage{amssymb,amsmath,eucal}

\begin{document}

\def\R {{\Bbb R }}
\def\C {{\Bbb C }}
\def\K{{\Bbb K}}
\def\H{{\Bbb H}}
\def\Z{{\Bbb Z}}

\newcommand{\arcsh}{\mathop{\rm arcsh}\nolimits}
\newcommand{\sh}{\mathop{\rm sh}\nolimits}

\def\const{{\rm const}}
\def\B{{\rm B}}

\def\OO{{\rm O}}
\def\SO{{\rm SO}}
\def\GL{{\rm GL}}
\def\SL{{\rm SL}}
\def\SU{{\rm SU}}
\def\Sp{{\rm Sp}}
\def\SOS{\SO^*}
\def\U{{\rm U}}

\def\ov{\overline}
\def\phi{\varphi}
\def\epsilon{\varepsilon}
\def\kappa{\varkappa}
\def\le{\leqslant}
\def\ge{\geqslant}

\def\densitybc{\Bigl|\frac
     {\Gamma(b+is)\Gamma(c+is)}{\Gamma(2is)}\Bigr|^2}

\def\densityabc{\Bigl|\frac
     {\Gamma(a+is)\Gamma(b+is)\Gamma(c+is)}{\Gamma(2is)}\Bigr|^2}

\def\densityx{x^{b+c-1}(1+x)^{b-c}}

\def\F{\,\,{}_2F_1}
\def\FF{\,\,{}_3F_2}

\def\FBC{\F(b+is, b-is; b+c; -x)}

\def\intt{\int_0^\infty}

\def\inttt{\int\limits_0^\infty}

\renewcommand{\Re}{\mathop{\rm Re}\nolimits}
\renewcommand{\Im}{\mathop{\rm Im}\nolimits}

\newcounter{sec}
 \renewcommand{\theequation}{\arabic{sec}.\arabic{equation}}

\newcounter{fact}
\def\fact{\addtocounter{fact}{1}{\scc \arabic{fact}}}

\newcounter{punkt}
\def\punkt{\addtocounter{punkt}{1}{\bff \arabic{punkt}}}

\begin{center}

{\Large\bf  Index hypergeometric transform and imitation
of analysis of Berezin kernels on hyperbolic spaces}

\bigskip

      {\sc Neretin Yu.A.}\footnote{ Supported by
     the grants RFBR 98-01-00303 and NWO 047--008--009}

\end{center}

\bigskip

{\small

The index hypergeometric transform (it is also called
 the Olevsky transform and the Jacobi transform)
generalizes the spherical transforms in  $L^2$
on rank 1 symmetric spaces (i.e. real, complex and quaternionic Lobachevsky
spaces). The purpose of the present work is to obtain
properties of the index hypergeometric transform
imitating the analysis of Berezin kernels on
rank 1 symmetric spaces.

We also discuss a problem of explicit construction of
a unitary operator identifying the space $L^2$
and a Berezin space. The problem is reduced to some integral
expression (the $\Lambda$-function)
that apparently can't be expressed in a finite form
in terms of standard special functions
(only for some special values of parameters this expression
can be reduced to the so-called Volterra type special functions).
We investigate some properties of this expression.
We show that for some series of symmetric spaces of large rank
the operator  of a unitary equivalence
can be expressed through the determinant of some matrix
consisting of the $\Lambda$-functions.
 }

\bigskip

{\bf 0.1. Mehler--Weyl--Titchmarsh--Olevsky
index hypergeometric transform.}
Fix $b,c>0$.  Consider the integral transform
$J_{b,c}$ defined dy the formula
\begin{multline}
 g(s)=J_{b,c}f(s)=\\=
 \frac 1 {\Gamma(b+c)}\intt f(x) \FBC\densityx dx
,\end{multline}
where $\F(\cdot)$ is the Gauss hypergeometric function.
The inverse operator  $J_{b,c}^{-1}$
   is given by
\begin{multline}
f(x)=  J_{b,c}^{-1} g(x)= \\ =
\frac 1{\pi\Gamma(b+c)}  \intt  g(s) \FBC
\densitybc ds
.\end{multline}

This pair of the inverse integral transforms
was discovered by Weyl
\cite{Wey}
in 1910. \footnote{The partial case $b=1/2, c=1$
(Mehler--Fock transform) was discovered by Mehler
\cite{Meh} in 1881.}
 Weyl's result didn't attract serious
interest and again these transforms appeared
in the works of Titchmarsh
 \cite{Tit}
in 1946
and M.N.Olevsky  \cite{Ole} в 1949г.
Spherical transforms for all Riemannian
noncompact symmetric spaces of rank 1 have the form
(0.1) for some special values of the parameters $b.c$.
This was the main reason of the interest  in
the transform (0.1) in subsequent years.

 For a detailed information on the transform
 (0.1) see the papers of Koornwinder and Flensted-Jensen
 \cite{FK}, \cite{Koo1}--\cite{Koo3},
see also the books of Yakubovich, Luchko
 \cite{YaL},  Samko, Kilbas, Marichev
 \cite{SKM} and   Tables of Prudnikov, Brychkov, Marichev
 \cite{PBM}, vol.5. Koornwinder's survey
 \cite{Koo2} contains a large bibliography.

For various generalizations of the transform (0.1)
see, for instance, the works of Yakubovich, Luchko
\cite{YaL} (the Wimp transform with $G$-function),
Koelink, Stockman
 \cite{Koe} ($q$-analogs),
Heckman, Opdam
\cite{HO}, \cite{Hec}, \cite{Opd}
(multivariate analogues).

In  mathematical literature there exist the following terms
for the transform (0.1)

--- the index hypergeometric transform

--- the generalized Fourier transform

--- the Olevsky transform

--- the Jacobi transform\footnote{The term was introduced
by Koornwinder by an association with expansions in
Jacobi polynomials.}

--- the Fourier--Jacobi transform

\smallskip

{\bf 0.2. Contents.} The space $L^2$ on
a noncompact Riemannian symmetric space $G/K$
admits a natural deformation discovered
(in the case of Hermitian symmetric spaces) by Berezin \cite{Ber2},
see also \cite{NO}, \cite{OO}, \cite{vDH}, \cite{Hil}
\cite{Ner1}--\cite{Ner3} and a bibliography in \cite{Ner3}.
For small values of the deformation parameter
(or for large values of our parameter  $\theta$, see \S 3)
the Berezin representations of the group
 $G$
are equivalent to the natural representation of $G$ in $L^2(G/K)$.

The Plancherel formula for the Berezin deformation
for rank 1 spaces was obtained by van Dijk and Hille
in
 \cite{vDH} and in the general case by the author  \cite{Ner1}--\cite{Ner3}.

One of purposes of the present work is to apply these results
to the theory of the index hypergeometric transform.

The second purpose is to obtain an explicit construction
of a unitary intertwining operator from $L^2(G/K)$ to
a Berezin space for large values of the parameter
$\theta$.
We also apply the theory of the index hypergeometric transform
 for  investigation of this operator.

Our \S 1 contains preliminaries on the index hypergeometric transform.

 In \S 2 we calculate the images of the operators
$$f(x)\mapsto xf(x);\qquad\qquad f(x)\mapsto x(x+1)\frac d{dx} f(x)$$
under the index hypergeometric transform. We show that these operators
correspond to difference operators in the imaginary direction
 with respect to the variable $s$.
 For instance,
 the image of the multiplication by $x$  is
\begin{multline*}
Pg(s)=\frac{(b-is)(c-is)}  {(-2is)(1-2is)}
    (g(s+i)-g(s))+
\frac{(b+is)(c+is)}  {(2is)(1+2is)}
    (g(s-i)-g(s))
\end{multline*}

Although these results are very simple,
I couldn't find
the statement on $x(x+1)\frac d{dx}$
in a literature.
The  statement on $x$ is a very partial case of one recent
Cherednik theorem
 \cite{Che}.

In \S 3 we give a brief survey of analysis of Berezin kernels
on rank 1 symmetric spaces.
Below we imitate this section on the level of the
index hypergeometric transform.

In \S 4 we give preliminaries on the continuous dual Hahn polynomials.

 In \S 5 we obtain a family of nonstandard Plancherel formulas.
Precisely, we show that the index transform is a unitary operator
from some space of holomorphic functions in a disk (or in a
 half-plane) to the space $L^2$ on a half-line with respect to some measure.
This section is an imitation of the works
\cite{vDH}, \cite{Ner1} on a level of special functions of
one variable.

 In \S 6 we define the function
\begin{multline*}
\Lambda^a_{b,c}(x)=\frac 1{\pi\Gamma(b+c)} \inttt
\frac{\Gamma(a+is)\Gamma(b+is) \Gamma(b-is) \Gamma(c+is) \Gamma(c-is)}
      {\Gamma(2is)\Gamma(-2is)} \times\\ \times
\FBC\,
ds
\end{multline*}
Our purpose is to construct  a unitary intertwining operator from
$L^2$ to a Berezin space.
Apparently, the function
 $\Lambda$ can't be expressed explicitly by means of
standard special functions  (with exception of some special values of $b$, $c$,
when $\Lambda$ can be reduced to the rare Volterra type
special functions).
We also try to understand an answer to the following informal
question: is it natural to consider the $\Lambda$-function
as a "new" special function? We find a collection
of identities with  $\Lambda$.
In \S 7 we construct some orthogonal bases
in $L^2$ on the half-line by means of the  $\Lambda$-function.

In \S 8 we present (without a proof)
 an explicit  construction of  a unitary operator
from $L^2$ to the Berezin deformation of $L^2$ for the
symmetric spaces
$$\U(p,q)/\U(p)\times\U(q).$$
It turned out that this integral operator  is expressed
by a determinant  of a matrix consisting of
the $\Lambda$-functions.
This section explains the purposes of introduction of
the $\Lambda$-function.

\smallskip

{\bf 0.3. Notation.}  $\R_+$ denotes the  half-line $x\ge0$.

$$\F(a,b;c;x):=\sum_{n=0}^\infty \frac{(a)_n(b)_n}{(c)_n \,n!}x^n$$
is the Gauss hypergeometric function.

$${}_pF_q\left[\begin{array}{c}a_1,\dots a_p\\c_1,\dots,c_q\end{array};x
  \right]:= \sum_{n=0}^\infty
  \frac{(a_1)_n\dots (a_p)_n} {(c_1)_n\dots (c_q)_n\, n!} x^n$$
is the generalized hypergeometric function.

\bigskip

{\large\bf \S 1. Preliminaries}

\addtocounter{sec}{1}
\setcounter{equation}{0}

\nopagebreak

\bigskip

{\bf 1.1. Definition.}
Fix $b,c>0$. Let $f(x)$ be a function on the half-line $\R_+$.
Its {\it index hypergeometric transform}
$J_{b,c}f$ (we also use the notation $[\widehat f]_{b,c}$)
is given by
\begin{multline}
 g(s)=J_{b,c}f(s)=[\widehat f(s)]_{b,c} =\\=
 \frac 1 {\Gamma(b+c)}\intt f(x) \FBC\densityx dx
.\end{multline}
The inversion formula for the index hypergeometric transform is
\begin{multline}
f(x)=  J_{b,c}^{-1} g(x)= \\ =
\frac 1{\pi\Gamma(b+c)}  \intt  g(s) \FBC
\densitybc ds
.\end{multline}
This is equivalent to the following statement:

  $J_{b,c}$ is a unitary operator
$$L^2\Bigl(\R_+, \densityx dx\Bigr)\to
  L^2\Bigl(\R_+, \densitybc dx\Bigr).$$

This property ({\it the Plancherel formula})
in the explicit form is
$$
\inttt f_1(x)\overline{f_2(x)}\densityx dx=
\frac 1\pi
\inttt [\widehat f_1(s)]_{b,c} \overline{[\widehat f_2(s)]}_{b,c}
 \densitybc ds
.$$
Behavior of the operator $J_{b,c}$ in different
functional spaces are known in details, see the survey
 \cite{Koo2}, and also \cite{YaL}.

\smallskip

{\sc Remark.} Weyl's (\cite{Wey}) proof
of the inversion formula
(1.2)
is relatively simple.
Consider the hypergeometric differential operator
$$D:=x(x+1)\frac {d^2}{dx^2}+\bigl[(c+b)+(2b+1)x\bigr] \frac {d}{dx}
+b^2$$
on $\R_+$.
 The functions $\FBC$ are the eigenfunctions of $D$.
There exists an explicit formula for the
resolution $R(\lambda)=(D-\lambda)^{-1}$  of the operator $D$
(this is equivalent to an explicit expression for the  Green function,
see  Myller-Lebedeff's works \cite{ML1}--\cite{ML2}).
If we know the resolution, then we can easily find the spectral
expansion by the formula
$$D=\frac 1{2\pi i}
 \int \lim\limits_{\epsilon\to 0}\bigl[ R(\lambda+i\epsilon)-
                R(\lambda-i\epsilon)\bigr]\,d\lambda
$$
A detailed  discussion of this reception
is contained in Dunford, Schwartz \cite{DS}.
For differential operators
(and in particular for the hypergeometric operator)
this method is present in  Titchmarsh's book \cite{Tit}.

For other proofs of the inversion formula see
 \cite{Ole}, \cite{Koo1}--\cite{Koo3},
\cite{SKM}, \cite{YaL}.

{\sc Remark.} Obviously,
$$[\widehat {Df}]=s^2 [\widehat {f}(s)]
$$
\smallskip

{\bf 1.2. Holomorphic extension to the strip.}

{\sc Lemma 1.1.} {\it
Let $f$ be integrable on $\R_+$ and
$$f(x)=o(x^{-\alpha-\epsilon}), \qquad\qquad x\to+\infty$$
where $\epsilon>0$. Then $[\widehat f(s)]_{b,c}$
is holomorphic in the strip}
$$|\Im s|<\alpha-b$$

{\sc Proof.}
This is a consequence  of the following asymptotics for the hypergeometric
function
(see \cite{HTF}, vol. 1, (2.3.2.9)) as $x\to+\infty$
$$
\FBC= \lambda_1 x^{-b+is} + \lambda_2 x^{-b-is}
+O(x^{-b+is-1}) +O(x^{-b-is-1}), $$
where $2s\not\in \Z$ and $\lambda_1$, $\lambda_2$ are some constants
(for  $2s\in \Z$ there appears the factor
$\ln x$ in the front of the leading term).

\smallskip

{\sc Remark.}
Assume a continuous function $f(x)$   satisfies
the condition
$$f(x)=o(x^{-b-\epsilon}), \qquad\qquad x\to+\infty.$$
Then
$f(x)$ is  an element of $L^2(\R_+)$ with respect to
the measure  $\densityx dx$.
We stress that this condition is also sufficient for existence
of analytic continuation to a narrow strip.

\smallskip

{\bf 1.3. Fourier transform.} For some special values
of the parameters
$b,c$ the function $\FBC$ is an elementary function. In particular
(see \cite{HTF}, vol. 1, (2.8.11--12)),
\begin{align}
\F(is,-is; \tfrac12;-x)=\cos(2s\arcsh\sqrt x),\\
\F(\tfrac12+is,\tfrac12-is; \tfrac32;-x)=    \frac{\sin(2s\arcsh\sqrt x)}
                                               {2s\sqrt x}
.\end{align}
Hence,
\begin{align*}
J_{0,1/2}f(s)&=\frac 1 {\Gamma(1/2)}
\intt f(x) \cos(2s\arcsh\sqrt x) x^{-1/2} (1+x)^{-1/2} dx, \\
J_{1/2,1}f(x)&=\frac 1 {\Gamma(3/2)} \intt f(x)
\frac{\sin(2s\arcsh\sqrt x)}
     {2s\sqrt x}
x^{1/2} (1+x)^{-1/2} dx
.\end{align*}
The substitution $x=\sh^2 y$ reduces these integrals to the form
\begin{align*}
J_{0,1/2}f(s)&=\frac 2 {\Gamma(1/2)}
 \intt f(\sh^2 y)\cos(2sy)dy,\\
J_{1/2,1}f(s)&=\frac 1 {\Gamma(3/2)}
\intt   f(\sh^2 y) \frac{\sh y\sin(2sy)}{2s}dy
,\end{align*}
i.e. the operator $J_{0,1/2}$ (resp.  $J_{1/2,1}$) differs from
the $\cos$-Fourier transform (resp. the $\sin$-Fourier transform)
by a nonessential variation of the notation.

Notice that the following functions also are elementary
$$
\F(k+is,k-is;\tfrac12+k+l;-x);\quad
\F(\tfrac12+k+is,\tfrac12+k-is;\tfrac32+k+l;-x)
.$$
for integer values of $k,l$. They can be obtained
from (1.3)--(1.4) by an application of suitable explicit
differential operators
(see \cite{HTF}, vol. 1, 2.8, formulas (20), (22), (24), (27)).

\smallskip

{\bf 1.4. Spherical transform.} Denote by  $\K$
the real numbers $\R$, the complex numbers $\C$,
or the quaternions $\H$.
Denote by $r$ the dimension of the field $\K$.
 By  $\K^n$ we denote the
$n$-dimensional space over  $\K$  with the standard
scalar product
$$\langle z,u\rangle=\sum z_j\overline u_j.$$
Denote by $\U(1,n;\K)$ the {\it pseudounitary group} over $\K$,
i.e. the group of all  $(1+n)\times(1+n)$ matrices
$\bigl(\begin{smallmatrix}a&b\\c&d\end{smallmatrix}\bigr)$ over $\K$,
satisfying the condition
$$
\begin{pmatrix}a&b\\c&d\end{pmatrix}
\begin{pmatrix}-1&0\\0&1\end{pmatrix}
\begin{pmatrix}a&b\\c&d\end{pmatrix}^*=
\begin{pmatrix}-1&0\\0&1\end{pmatrix}
.$$
The standard notations for  $\U(1,n;\K)$
in the cases $\K=\R,\C,\H$ are respectively:
$\OO(1,n)$, $\U(1,n)$, $\Sp(1,n)$.

By $\B_n(\K)$ we denote the unit ball
$\langle z,z\rangle <1$ in $\K^n$. Let $S^{rn-1}$ be the
unit sphere
$\langle z,z\rangle=1$. The group $\U(1,n;\K)$ acts on $\B_n(\K)$
by the linear fractional transformations
\begin{equation}
z\mapsto z^{[g]}:=(a+zc)^{-1}(b+zd)
.\end{equation}
The stabilizer $K$ of the point $0\in \B_n(\K)$
consists of the matrices
\begin{equation}
\begin{pmatrix}a&0\\0&d\end{pmatrix}
                     \qquad\qquad |a|=1,\quad d\in\U(n;\K)
.\end{equation}
Hence   $\B_n(\K)$  is a homogeneous space,
$$
\B_n(\K)=\U(1,n;\K)/\U(1;\K)\times\U(n;\K)
.$$

The Jacobian of  transformation
 (1.5) is
$$
J(g;z)=|a+zc|^{-r(1+n)}
.$$
The following formula can  easily be checked
$$
1-\langle z^{[g]},u^{[g]} \rangle
=(a+zc)^{-1}(1-\langle z,u\rangle)\overline{ (a+uc)}\,^{-1}
.$$
This implies that the $\U(1,n;\K)$-invariant measure
on $\B_n(\K)$
has the form
$$
dm(z)=(1-\langle z,z\rangle)^{-(n+1)r/2} dz
,$$
where $dz$ denotes the Lebesgue measure on    $\B_n(\K)$.

The group
$\U(1,n;\K)$
acts in  $L^2(\B_n(\K),  dm(z))$ by the substitutions
\begin{equation}
f(z)\mapsto f\bigl((a+zc)^{-1}(b+zd) \bigr)
.\end{equation}

Consider the problem of the decomposition of
 $L^2(\B_n(\K),  dm(z))$ into
a direct integral of irreducible representations of
the group
 $\U(1,n;\K)$.

Let $s\in\R$.
The representation $T_s$ of (spherical) {\it principal unitary series}
  of the group $\U(1,n;\K)$
acts in  $L^2(S^{rn-1})$ by
the operators
\begin{equation}
T_s
\begin{pmatrix}a&b\\c&d\end{pmatrix}
f(h)=f\bigl( (a+hc)^{-1}(b+hd)\bigr)|a+hc|^{-(n+1)r/2+1+is}
,\end{equation}
where $h\in S^{rn-1}$.

Consider the operator
 $A$ from the space
          $L^2(\B_n(\K),  dm(z))$
to the space of functions on $S^{rn-1}\times \R_+$
given by
\begin{equation}
Af(h,s)=\int_{\B_n(\K)}            f(z)
\frac{|1-\langle z,h \rangle|^{-(n+1)r/2+1+is}      }
     {|1-\langle z,z \rangle|^{(n+1)r/4+1/2+is/2} }   dz
.\end{equation}
If $f$ is transformed by (1.7), then
$G(h,s):=Af(h,s)$ is transformed by
\begin{equation}
G(h,s)\mapsto
G\bigl( (a+hc)^{-1}(b+hd),\,\,s\bigr)|a+hc|^{-(n+1)r/2+1+is}
.\end{equation}

Now we want to construct the measure
  $d\nu$ (it is called   the
{\it Plancherel measure}) on
$S^{rn-1}\times \R_+$ such that the operator
 $A$ is a unitary operator
$$L^2 (\B_n(\K),  dm(z))  \to L^2(S^{rn-1}\times \R_+, d\nu)$$
(this will solve the  problem of decomposition into a direct integral).
Obviously, the factor $|a+hc|$ (см.(1.8),(1.9))
is 1 if the matrix
$\bigl(\begin{smallmatrix}a&b\\c&d\end{smallmatrix}\bigr)$
has the form (1.6). Thus the required measure $d\nu$ is
invariant with respect to the group
 $\U(n;\K)$. Hence the measure $d\nu$ has the form
$$\sigma(s) \,ds \,dh.$$
To find the function $\sigma(s)$,
we consider all functions
$f$ on $\B_n(\K)$
depending only on the radius  $|z|$. It is convenient to introduce the
variable
$$x=\frac{|z|^2}{1-|z|^2}$$
and to assume $f=f(x)$. Then the corresponding function $G(h,s)$
 depends only in the variable $s$ and after a simple calculation
we obtain
\begin{equation}
G(s)=\const\cdot\intt f(x) \FBC\densityx dx
,\end{equation}
where
\begin{equation}
b=(n+1)r/4-1/2;\qquad c=(n-1)r/4+1/2
.\end{equation}
The integral transform (1.11) defined on
the space of $\U(n,\K)$-invariant functions
on $\B_n(\K)$
 is called the {\it spherical transform}.
We observe that the spherical transforms
are special cases of the index hypergeometric transform.
The inversion formula (1.2) gives the following density
for the Plancherel measure
\begin{equation}
\sigma(s)=\densitybc
\end{equation}

\bigskip

       {\large\bf \S 2. Correspondence of some operators}

\nopagebreak

\bigskip

\addtocounter{sec}{1}
\setcounter{equation}{0}

{\bf 2.1. Operator of multiplication by $x$.}

{\sc  Theorem 2.1.} {\it Assume that a function  $f$ defined  on $\R^+$
is continuous and satisfies the condition
\begin{equation}
f(x)=o(x^{-b-1-\epsilon});\qquad x\to +\infty
.\end{equation}
Then
\begin{equation}
[\widehat{xf(x)}]_{b,c}= P[\widehat{f(x)}]_{b,c}
,\end{equation}
where the difference operator  $Pg$ is given by}
\begin{multline*}
Pg(s)=\frac{(b-is)(c-is)}  {(-2is)(1-2is)}
    (g(s+i)-g(s))+
\frac{(b+is)(c+is)}  {(2is)(1+2is)}
    (g(s-i)-g(s))
.\end{multline*}

{\sc Proof.} By condition (2.1), there exists a holomorphic continuation
of the function   $\widehat{f}$ to the strip
$$|\Im s|<1+\epsilon.$$
The identity (2.2) is equivalent to the following identity
for hypergeometric functions
\begin{multline*}
x\FBC=\\
=\frac{(b-is)(c-is)}  {(-2is)(1-2is)}
  \F(b-1+is,b+1-is;b+c;-x)-\\
-\Bigl[\frac{(b-is)(c-is)}  {(-2is)(1-2is)}+
     \frac{(b+is)(c+is)}  {(2is)(1+2is)}\Bigr]\FBC+\\
+\frac{(b+is)(c+is)}  {(2is)(1+2is)}
 \F(b+1+is,b-1-is;b+c;-x)
.\end{multline*}
The last identity is equivalent to the following relation
between associated hypergeometric functions%
\footnote{For any triple of hypergeometric functions
$\F(p,q;r;y)$, $\F(p+k_1,q+l_1;r+m_1;y)$,
$\F(p+k_2,q+l_2;r+m_2;y)$ with integer $k_j,l_j,m_j$
there exists a linear relation with  coefficient,which are
  polynomial in $y$s.}
\begin{multline}
-y \F(p,q;r;y)=\\
=\frac{q(r-p)}{(q-p)(1+q-p)}\F(p-1,q+1;r;y)-\\
-\Bigl[\frac{q(r-p)}{(q-p)(1+q-p)}+
 \frac{p(r-q)}{(p-q)(1+p-q)}\Bigr]\F(p,q;r;y)+\\
+\frac{p(r-q)}{(p-q)(1+p-q)}\F(p+1,q-1;r;y)
.\end{multline}
The last identity is absent in standard tables
but it can easily  be checked by  equating of the coefficients of
the Taylor
expansion in
 $y$.
Indeed,
\begin{multline*}
\F(p-1,q+1;r;y)-\F(p,q;r;y)
=\sum\left( \frac{(p-1)_k(q+1)_k}{(r)_k k!}
-\frac{(p)_k(q)_k}{(r)_k k!}\right)=\\
=\sum  \frac{(p)_k(q)_k}{(r)_k k!}\left(\frac{(p-1)(q+k)}{(p+k-1)q}-1\right)=
\sum  \frac{(p)_k(q)_k}{(r)_k k!} \cdot \frac{-k(1+q-p)}{q(p+k-1)}
.\end{multline*}
Converting
$\F(p+1,q-1;r;y)-\F(p,q;r;y)$,
we reduce the right part of (2.3) to the form
$$\sum  \frac{(p)_k(q)_k}{(r)_k k!}
    \left( \frac{-k(r-p)}{(q-p)(p+k-1)} -
\frac{-k(r-q)}{(p-q)(q+k-1)}\right).$$
The last expression equals the left part of
  (2.3)

\smallskip

{\bf 2.2.  Differentiation.}
{\sc Theorem 2.2.} {\it Let $f$, $f'$ be continuous
and satisfying the decreasing conditions  {\rm (2.2)}.
Then
\begin{equation}
[\widehat{x(x+1)\frac{d}{dx}f}]_{b,c}=H[\widehat f]_{b,c}
,\end{equation}
where the difference operator $H$ is given by }

\begin{multline*}
Hg(s)
=
\frac{(b-is)(b+1-is)(c-is)}  {(-2is)(1-2is)}
         (g(s+i)-g(s))+\\
+
\frac{(b+is)(b+1+is)(c+is)}  {(+2is)(1+2is)}
         (g(s-i)-g(s))-(b+c)g(s)
.\end{multline*}

{\sc Proof.}
\begin{align*}
\intt\Bigl\{ x(x+1)\frac{d}{dx}f(x)\Bigr\}\FBC\densityx dx=\\
=-\intt f(x)\frac{d}{dx}\Bigl\{x^{b+c}(1+x)^{b-c+1}\FBC\Bigr\}dx=\\
=-\intt f(x) \Bigl\{x(x+1)\frac{d}{dx}\FBC+\\
+\bigl((2b+1)x+(b+c)\bigr)\FBC\Bigr\}\densityx dx
.\end{align*}
Now the statement is reduced to the identity
\begin{align*}
&-\Bigl\{x(x+1)\frac{d}{dx}+(2b+1)x+(b+c)\Bigr\}\FBC=\\
=&
\frac{(b-is)(b+1-is)(c-is)}  {(-2is)(1-2is)}
                                \F(b-1+is,b+1-is;b+c;-x)
-\\&- \left(
\frac{(b-is)(b+1-is)(c-is)}  {(-2is)(1-2is)}
+\frac{(b+is)(b+1+is)(c+is)}  {(+2is)(1+2is)}
-b-c\right)\times\\
&\qquad\qquad\qquad\qquad\qquad\qquad\qquad\qquad\times \FBC+
\\&
+\frac{(b+is)(b+1+is)(c+is)}  {(+2is)(1+2is)}
\F(b+1+is,b-1-is;b+c;-x)
.\end{align*}
This is an relation between associated hypergeometric functions
\begin{multline*}
 -\Bigl\{x(x+1)\frac{d}{dx}+(p+q+1)x+r\Bigr\}\F(p,q;r;-x)=\\
=(x^2+x)\frac{pq}r\F(p+1,q+1;r+1;-x)- \bigl((p+q+1)x+r\bigr)\F(p,q;r;-x)=\\
=\frac{q(q+1)(r-p)}{(1+q-p)(p-q)}
\bigl( \F(p-1,q+1;r;-x)-\F(p,q;r;-x) \bigr)+
\\+
\frac{p(p+1)(r-q)}{(1+p-q)(q-p)}
\bigl( \F(p+1,q-1;r;-x)-\F(p,q;r;-x) \bigr)
\end{multline*}
and can easily be checked by equating of the Taylor coefficients at
 $x^k$.

\bigskip

{\large\bf \S3. Berezin kernels}

\nopagebreak

\bigskip

\addtocounter{sec}{1}
\setcounter{equation}{0}

 The space $L^2$ on the ball $\B_n(\K)$
(the ball $\B_n(\K)$ is the symmetric space
$\U(1,n;\K)/\U(1,\K)\times\U(n,\K)$)
admits a natural deformation that will be described in this section
(for more details see  \cite{vDH} for rank 1 symmetric spaces
and \cite{Ner3}
for symmetric spaces of arbitrary rank).
The main purpose of the present work
is an imitation of this deformation on a level of special functions.

\smallskip

{\bf 3.1. Positive definite kernels.}
Recall the definition of a {\it  positive definite kernel}
(for detailed discussion and references see \cite{Ner3}).
Let $X$ be a set. A function $L(x,y)$  on $X\times X$
is called {\it a positive definite kernel} if
$L(y,x)=\overline{L(x,y)}$ and for all $x_1,\dots, x_k\in X$
the matrix
$$\det\begin{pmatrix}
L(x_1,x_1)&\dots& L(x_1,x_k)\\
\vdots&\ddots&\vdots \\
L(x_k,x_1)&\dots& L(x_k,x_k)
\end{pmatrix}$$
is positive semidefinite.

Let $L$ be a positive definite kernel on a set $X$.
Then there exists a Hilbert space $H=H[L]$
and a system of vectors $v_x\in H$, where $x$ ranges $X$,
such that

 1. the linear span of the vectors $v_x$ is dense in $H[L]$

 2. $\langle v_x, v_y \rangle =L(y,x)$.

\smallskip

For each vector
$h\in H[L]$ we assign
the function $f_h(x)$ on $X$ by
$$f_h(x)=\langle h,v_x \rangle$$
Thus the space
 $H$ is realized as some space of functions on
$X$. We denote this function  space by
$H^\circ[L]$. Obviously, the vectors $v_a\in H[L]$
correspond to the functions
$$\phi_a(x)=L(x,a)$$
and for each $f\in H^\circ[L]$ we have the following equality
 (the {\it reproducing property})
\begin{equation}
\langle f, \phi_a\rangle_{H^\circ[L]}=f(a)
.\end{equation}
This identity gives some (constructive, but not completely)
description of the scalar product in
$H^\circ[L]$.
It is said also that $L$ is a {\it reproducing kernel}
of the space
$H^\circ[L]$.

\smallskip

{\bf 3.2. Berezin kernels.}
The {\it Berezin kernel}
 $L_\theta(z,u)$ on $\B_n(\K)$ is defined by
\begin{equation}
L_\theta(z,u)=|1- \langle z,u \rangle|^{-\theta}
.\end{equation}
In the cases $\K=\R,\C$ the kernel $L_\theta(z,u)$
is positive definite
for all $\theta\ge0$, and in the case $\K=\H$
the condition of positive definiteness is $\theta\ge 2$.
The group $\U(1,n;\K)$ acts in $H^\circ(L_\theta)$
by the unitary operators
$$T_\theta
\begin{pmatrix}a&b\\c&d\end{pmatrix}
f(z)=f\bigl((a+zc)^{-1}(b+zd)\bigr) |a+zс|^{-\theta}
.$$

\smallskip

{\bf 3.3. Another version of the definition.}
Consider the kernel on $\B_n(\K)$ given by
$$M_\theta(z,u)=\frac{(1- \langle z,\,z\rangle)^{\theta/2}
   (1-\langle u,\,u\rangle)^{\theta/2}}
   {|1-\langle z,\,u\rangle|^{\theta}}
.$$
Consider the spaces
$H(M_\theta)$, $H^\circ(M_\theta)$ associated with this kernel.
The group $\U(1,n;\K)$ acts in the space $H^\circ[M_\theta]$
by substitutions by the formula
\begin{equation}
R_\theta
\begin{pmatrix}a&b\\c&d\end{pmatrix}
f(z)=f\bigl((a+zc)^{-1}(b+zd)\bigr)
.\end{equation}
It is easily shown that the constructions of 3.2 and  3.3
are equivalent.
The canonical unitary  $\U(1,n;\K)$-intertwining operator
$H^\circ[L_\theta] \to H^\circ[M_\theta]$
is the operator of the multiplication by the function
$(1-\langle z,\,z \rangle)^{\theta/2}$.

\smallskip

{\bf 3.4. More material description of the space
 $H^\circ(L_\theta)$.}
The cases $\K=\R,\C,\H$ are similar and we will
consider the case
$\K=\R$ as a basic example. Let us give a more constructive
description of the space
$H^\circ(L_\theta)$ in this case
(with an additional restriction $\theta>n$).

For this aim we consider the space $V_\theta$
of holomorphic functions
in the ball $\B_n(\C)$, we  equip this space with the scalar
product
$$\langle f,g\rangle_{V_\theta} =\frac{\Gamma(\theta)}{\pi^n\Gamma(\theta-n)}
\int_{\B_n(\C)} f(z)\overline{g(z)}(1-\langle z,z\rangle)^{\theta-n-1} dz
.$$

The universal covering group
of the group $\U(1,n)=\U(1,n;\C)$ acts in the space $V_\theta$
by the unitary operators
$${\cal U}_\theta \begin{pmatrix} a&b\\c&d\end{pmatrix}
f(z)= f\bigl((a+zc)^{-1}(b+zd)\bigr) (a+zc)^{-\theta}.
$$

{\sc Remark.} For an integer $\theta$
the family of the operators ${\cal U}_\theta$
defines a linear representation of the group $\U(1,n)$.
 If $\theta$ is not integer,
then
$$(a+zc)^{-\theta}=
       (1+zca^{-1})^{-\theta}e^{-\theta(\ln a+2\pi k i)}.
$$
It can  easily be shown, that $\langle ca^{-1}, ca^{-1}\rangle<1$,
hence $(1+zca^{-1})^{-\theta}$ is canonically defined
for $\langle z,z\rangle<1$. Thus the operator
${\cal U}_\theta
\left( \begin{smallmatrix} a&b\\c&d\end{smallmatrix}\right)$
is defined up to a factor $e^{2\pi ki\theta}$.
Hence for a noninteger $\theta$
we obtain a projective representation
of the group $\U(1,n)$ or (this is one and the same)
a linear representation of its universal covering group.

\smallskip

A simple calculation (with the Dirichlet integral, see,
for instance, \cite{AAR}, 1.8)
shows that the system of functions
$z_1^{k_1}\dots z_n^{k_n}$ forms  an orthogonal basis in $V_\theta$
and
\begin{equation}
\| z_1^{k_1}\dots z_n^{k_n} \|^2
  =\frac{ k_1!\dots k_n!}{ (\theta)_{k_1+\dots+ k_n}}
.\end{equation}
Consider the function
\begin{equation}
\phi_a(z):=\bigl(1-\sum z_j \overline a_j\bigr)^{-\theta}=
\sum_{k_1,\dots,k_n}
\frac { (\theta)_{k_1+\dots k_n}
\overline a_1^{k_1}\dots \overline a_n^{k_n} } { k_1!\dots k_n!}
z_1^{k_1}\dots z_n^{k_n}
.\end{equation}
By equations (3.4)--(3.5), for any  function $f\in V_\theta$
$$
\langle f, \phi_a\rangle_{V_\theta}=f(a)
.$$
Thus the space $V_\theta$ is defined by the positive definite kernel
$$
K(z,u)=\bigl(1-\langle z,\,u\rangle \bigr)^{-\theta}
.$$
A function holomorphic in the ball $\B_n(\C)$
is uniquely determined by its restriction
to         $\B_n(\R)$,
hence we can consider the space $V_\theta$ as a space
of functions on
$\B_n(\R)$.
This is strictly the space defined by the Berezin kernel
 $L_\theta$ on $\B_n(\R)$.

Emphasize that in the space
 $H^\circ[L_\theta]$
on $\B_n(\R)$ we have an action of the group $\U(1,n)$,
which is larger than
$\OO(1,n)$\footnote{ For the complex ball $\B_n(\C)$,
 the  analogous overgroup
is  $\U(1,n)\times\U(1,n)$.
For the quaternionic ball the overgroup is
$\U(2,2n)$}; for more details see \cite{Ner3}.

\smallskip

{\bf 3.5.  Plancherel formula.} Emphasize, that the group
$\U(1,n;\K)$
acts in the space $L^2(\B_n(\K))$ and in all spaces
$H^\circ(M_\theta)$ by the same formula (1.7), (3.3).
Consider the following operator, which coincides with
 operator (1.9) up to a functional factor
\begin{multline}
R_\theta f(h,s)=\frac 1{\Gamma(b+c)|\Gamma(\theta-b+is)|^2}
           \int\limits_{\B_n(\K)}            f(z)
\frac{|1-\langle z,h \rangle|^{-(n+1)r/2+1+is}      }
     {|1-\langle z,z \rangle|^{(n+1)r/4+1/2+is/2} }   dz
.\end{multline}
As was shown by van Dijk and Hille \cite{vDH},
 for $\theta>b$ the operator $R_\theta$ is a unitary operator
$$H^\circ(M_\theta)\to L^2(S^{rn-1}\times \R_+, \tau(s)\,dh\,ds),$$
where $dh$ is a Lebesgue measure on the sphere and
\begin{equation}
\tau(s) =\frac 1{\Gamma(\theta)\Gamma(\theta-b+c)}
   \left|
  \frac{\Gamma(\theta-b+is)\Gamma(b+is)\Gamma(c+is)}{\Gamma(2is)}
          \right|^2
.\end{equation}

 {\sc Remark.} If we remove the expression
 $|\Gamma^2(\theta-b+is)|^2$ from the denominator in (3.6),
then in (3.7) the same expression will relocated
from the numerator to the denominator.
By some reasons, the normalization (3.6)
is more convenient.

\smallskip

Below in \S 5 we obtain Theorem 5.3 imitating the Plancherel formula (3.7)
 on the level of the index hypergeometric transform,
the  van Dijk--Hille theorem is a special case of Theorem 5.3.

\smallskip

{\bf 3.6.  Radial parts of Berezin kernels.}
Again let us for definiteness consider the case
$\K=\R$. Consider the space $V_\theta$
of holomorphic functions defined in 3.4,
and its subspace $V_\theta^{\OO(n)}$
({\it the radial Berezin space})
consisting of $\OO(n)$-invariant functions.

Obviously, the elements of the space $V_\theta^{\OO(n)}$ have the form
$$
g(z_1^2+\dots+z_n^2)=\sum c_p (z_1^2+\dots+z_n^2)^p .
$$
Hence we can consider $V_\theta^{\OO(n)}$ as a space of functions depending
on one variable
$$u=z_1^2+\dots+z_n^2$$
 lying in the disk
$|u|<1$.

\smallskip

{\sc Lemma 3.1.} {\it  Vectors $u^p$ form an orthogonal basis
in
the space $V_\theta^{\OO(n)}$ and}
\begin{equation}
\|u^p\|^2=\frac{p!\, (n/2)_p}
               {(\theta/2)_p (\theta/2+1/2)_p}
.\end{equation}

\smallskip

{\sc Proof.}
Obviously, the vectors $u^p$ are pairwise orthogonal,
hence it is sufficient to calculate $\|u^p\|^2$.
By orthogonality (3.4),
we obtain

\begin{multline*}
\|u^p\|^2=\|(z_1^2+\dots+z_n^2)^p\|^2 =
\sum\limits_{k_1\ge 0,\dots,k_n\ge 0, \sum k_j=p}
\left(\frac{p!}{k_1!\dots k_n!}\right)^2 \|z_1^{2k_1}\dots z_n^{2k_n}\|^2
 =\\=
\sum\limits_{k_1\ge 0,\dots,k_n\ge 0, \sum k_j=p}
\left(\frac{p!}{k_1!\dots k_n!}\right)^2
\frac{2k_1!\dots 2k_n!}{(\theta)_{2p}}
.\end{multline*}
It is sufficient to evaluate
\begin{equation}
\sum\limits_{k_1\ge 0,\dots,k_n\ge 0, \sum k_j=p}
\frac{2k_1!\dots 2k_n!}{(k_1!\dots k_n!)^2}
.\end{equation}
Denote a summand of this sum by $A_{k_1,\dots, k_n}$.
Consider the generating function
\begin{multline*}
h(y_1,\dots y_n):=
\sum\limits_{k_1\ge 0,\dots,k_n\ge 0}
A_{k_1,\dots, k_n}y_1^{k_1}\dots y_n^{k_n} =\\=
\prod_{m=1}^n \sum_{k=0}^\infty \frac{2k!}{k!\,k!}y_m^k=
\prod_{m=1}^n (1-4y_m)^{-1/2}
.\end{multline*}
Expression (3.9) is the coefficient at $y^p$
in the Taylor expansion of
$$h(y,\dots,y)=(1-4y)^{-n/2}$$
and this gives the required statement.

\smallskip

Spaces that imitate $V_\theta^{\OO(n)}$
are defined below in 5.1.

\smallskip

{\bf 3.7. {\O}rsted problem.}
As before, let
 $b,c$ be defined by (1.12).
 For $\theta>b$ the representation
$\U(1,n;\K)$ in $H^\circ[L_\theta]$ is equivalent to the
representation
in $L^2(\B_n(\K))$. There arises the following question:

{\it Is it possible to write explicitly a unitary intertwining operator}
$$S:\,\,L^2(\B_n(\K)) \to  H^\circ[L_\theta]\qquad \mbox{\Large ?}$$

It is easy to obtain a general form of this operator.
The operator  $S$ is a product of three operators
$$S=\const\cdot
R_\theta^{-1} M A,$$
where $A$ is given by
(1.9), $R_\theta$ is defined by (3.6), and
$$M:L^2\left(\R_+,\densitybc\right)
   \to L^2\left(\R_+,
\Bigl|\frac {\Gamma(\theta-b+is)\Gamma(b+is)\Gamma(c+is)}
            {\Gamma(2is)}\Bigr|^2
\right)$$
is the operator of division by a function
$\psi(s)$ satisfying the condition
$$|\psi(s)|=|\Gamma(a-b+is)|.$$

{\sc Lemma 3.2.} {\it The operator $S$ is given by
$$Sf(z)=\int_{\B_n(\K)}
\Lambda\left(\frac{|z-u|^2}{(1-|z|^2)(1-|u|^2)}\right) f(u)
\,dm(u)
,$$
where $dm$ is the  invariant measure on the ball and}
\begin{equation}
\Lambda(x)=\const\cdot\intt\overline{\psi(s)}\densitybc\FBC \,ds
.\end{equation}

{\sc Proof.} Consider the function $1$ from
$V_\theta$. Its preimage with respect to the operator $S$
can be easily evaluated, it is given by (3.10),
 where $x=|z|^2/(1-|z|^2)$.
By the reproducing property,
for any
 $g\in V_\theta$ we have
$\langle g,1\rangle=g(0)$.  Since the operator $S$ is unitary,
 any function $f\in L^2$
satisfies
$$
Sf(0)=\int_{\B_n(\K)} \Lambda\Bigl(\frac{|z|^2}{1-|z|^2)}\Bigr)
                       dm(z)
.$$
Now the kernel of the operator $S$ can be reconstructed by
invariance arguments.

\smallskip

Where arises a question, is it possible to find
 $\psi(x)$ such that the integral (3.10)
can be explicitly evaluated.  For several years I tried
to solve this problem and now I think that this is impossible.
It seems that the simplest variant
$$\psi(s)=\Gamma(\theta-b-is)$$
is best.

In \S 6 I am trying to understand, is it natural
to consider  the integral (3.10) obtained in this way
as a 'new'   special function.

I am inclined to think that the answer is affirmative by the following
 'metaphysical' reason. To be definite, assume $\K=\R$.
Then the natural symmetry group   $\U(1,n)$
of the space $V_\theta$ is larger than the symmetry group
$\OO(1,n)$ of  the space $L^2(\B_n(\R))$.
If we identify $L^2$ with $V_\theta$, then we force the group
$\U(1,n)$ to act in $L^2(\B_n(\R))$.
It is natural to think that enlarging of a symmetry group
must imply a nontrivial analysis on the level of special functions.

There is also another fact that seems pleasant for me.
For 3 series of semisimple groups
$\OO(n,\C)$, $\Sp(n,\C)$,
$\U(p,q)$
the corresponding  $\Lambda$-function can be expressed
as a determinant consisting of $\Lambda$-functions (3.10)
related to rank 1 symmetric spaces.

\bigskip

{\large\bf \S4. Hahn polynomials: preliminaries.}

\nopagebreak

\bigskip

\addtocounter{sec}{1}
\setcounter{equation}{0}

{\bf 4.1. Definition.} Пусть $a,b,c>0$.
The {\it continuous dual Hahn polynomials}
(see \cite{AAR}, \cite{VK}) are given by
$$S_n(s^2;a,b,c)=(a+b)_n(a+c)_n
\FF\left[
  \begin{array}{c}-n, a+ is , a-is\\ a+b,a+c
   \end{array};1\right]
.$$
Obviously, the expression $S_n(s^2;a,b,c)$ doesn't change
if we if we transpose the parameters $b,c$.
The Kummer formula (see \cite{AAR}, Corollary 3.3.5.)
$$
\FF\left[
  \begin{array}{c}\alpha, \beta, \gamma\\ \delta, \epsilon
   \end{array};1\right]=
\frac{\Gamma(\epsilon)\Gamma(\delta+\epsilon-\alpha-\beta-\gamma)}
     {\Gamma(\epsilon-\alpha)\Gamma(\delta+\epsilon-\beta-\gamma)}
\FF\left[
  \begin{array}{c}\alpha,\delta- \beta,\delta- \gamma\\ \delta,
 \delta+\epsilon-\beta-\gamma
   \end{array};1\right]
$$
implies that $S_n(s^2;a,b,c)$ is invariant with respect to
arbitrary permutations       of
$a,b,c$.

\smallskip

{\bf 4.2. Orthogonality relations.}  The Hahn polynomials
form an orthogonal basis in $L^2$ on half-line
with respect to the weight
$$
\densityabc$$
and
\begin{multline*}
\frac 1{\pi}\intt
\densityabc S_n(s^2;a,b,c)S_m(s^2;a,b,c)
\,dx=\\=
\Gamma(a+b+n)\Gamma(a+c+n)\Gamma(b+c+n)\, n! \delta_{m,n}
.\end{multline*}

{\bf 4.3. Difference equations.}
Consider the difference operator
\begin{equation}
{\cal L}y(s)=B(s)y(s+i)-(B(s)+D(s))y(s)+D(s)y(s-i)
,\end{equation}
where
\begin{align*}
B(s)=
 \frac{(a-is)(b-is)(c-is)}{(-2is)(1-2is)},\\
D(s)=
 \frac{(a+is)(b+is)(c+is)}{(+2is)(1+2is)}
.\end{align*}

The Hahn polynomials are
the eigenfunctions of this operator
\begin{equation}{\cal L} S_n(s^2;a,b,c)= n\cdot S_n(s^2;a,b,c)
.\end{equation}

{\bf 4.4. Index hypergeometric transform and Hahn polynomials.}

{\sc Lemma 4.1.} {\it The image of the function $(1+x)^{-a-b}$
under the transform $J_{b,c}$ is}
$$\frac{\Gamma(a-is)\Gamma(a+is)}
  {\Gamma(a+b)\Gamma(a+c)}
$$

{\sc Proof.} The statement is reduced to the table integral
\cite {GR}, 7.51.10.

\smallskip

{\sc Lemma 4.2.} {\it The image of the function
$$\left(\frac x{x+1}\right)^n (1+x)^{-a-b}$$
with respect to the transform
 $J_{b,c}$ is}
\begin{multline*}
\frac{\Gamma(a-is)\Gamma(a+is)}
  {\Gamma(a+b)\Gamma(a+c)}
\FF\left[
  \begin{array}{c}-n, a+ is , a-is\\ a+b,a+c
   \end{array};1\right]
=\\            =
\frac{|\Gamma(a+is)|^2}{\Gamma(a+b+n)\Gamma(a+c+n)}
S_n(s^2;a,b,c)
.\end{multline*}

{\sc Proof.} We must evaluate the image of the function
\begin{multline*}
\left(\frac x{x+1}\right)^n (1+x)^{-a-b}=
\left(\frac {1+x-1}{x+1}\right)^n (1+x)^{-a-b}=\\=
\sum_{k=0}^n C_n^k(-1)^{n-k}(1+x)^{-a-c-n+k}
.\end{multline*}
But the images of the functions $(1+x)^{-a-c-n+k}$
were evaluated above.

\bigskip

{\large\bf \S 5. Nonstandard Plancherel formulas}

\nopagebreak

\bigskip

\addtocounter{sec}{1}
\setcounter{equation}{0}

{\bf 5.1. Spaces $H^a_{b,c}$.} Fix positive numbers
$a,b,c$ satisfying the conditions
\begin{equation}
a>b,\quad a>c,\quad 2a>1
.\end{equation}
Consider the space $W^a$ of holomorphic functions
in the disk
$|z|<1$ satisfying the condition
$$
\iint\limits_{|z|<1}\bigl|f(z)\bigr|^2(1-|z|^2)^{2a-2} dz\,<\infty
,$$
where
$dz$ denotes the Lebesgue measure in the disk.

Let us define the scalar product in
 $W^a$ by the formula
$$\langle f ,g\rangle=\frac 1{\pi\Gamma(2a-1)}
  \iint\limits_{|z|<1} f(z)\overline{g(z)}  (1-|z|^2)^{2a-2}
 \F(a-b, a-c; 2a-1; 1-|z|^2)\,dz.$$
We denote by   $W_{b,c}^a$ the Hilbert space obtained in this way.

\smallskip

{\sc Remark.} The  spaces $W^a_{b,c}$  as linear spaces
don't depend on $b,c$ and coincide with $W^a$.  But the scalar
products in these spaces are  different.

\smallskip

{\sc Lemma 5.1.} {\it The functions $z^k$ form an orthogonal
basis in $W_{b,c}^a$ and}
\begin{equation}
\langle z^k,\,z^k\rangle_{W_{b,c}^a}
 =\frac{k!\, \Gamma(b+c+k)}{\Gamma(a+b+k)\Gamma(a+c+k)}
.\end{equation}

{\sc Proof.} Orthogonality of the functions $z^k$ is obvious.
Let us evaluate
\begin{align*}
&\langle z^k,\,z^k\rangle=\\&=
\frac 1{\pi\Gamma(2a-1)}
  \iint\limits_{|z|<1} |z|^{2k}  (1-|z|^2)^{2a-2}
 \F(a-b, a-c; 2a-1; 1-|z|^2)\,dz
=\\
&=
\frac 2{\Gamma(2a-1)}
\int_0^1 r^{2k+1}(1-r^2)^{2a-2}
  \F(a-b, a-c; 2a-1; 1-r^2)dr =\\
&=
\frac 1{\Gamma(2a-1)}
\int_0^1 y^{k}(1-y)^{2a-2}
  \F(a-b, a-c; 2a-1; 1-y)dy =\\
&=
\frac 1{\Gamma(2a-1)}
\int_0^1 (1-v)^{k}v^{2a-2}
  \F(a-b, a-c; 2a-1; v)dv
\end{align*}
Then all is reduced to the table integral
 \cite{GR}, 7.512.4.

\smallskip

{\sc Remark.} We observe that the radial Berezin spaces
$V_\theta^{\OO(n)}$
defined in 3.6 are special cases of the spaces
$W_{b,c}^a$:
$$b=n/4-1/4;\qquad c=n/4+1/4;\qquad a=\theta-n/4+1/4.$$

\smallskip

{\bf 5.2. Reproducing kernels of the spaces $W_{b,c}^a$.}

{\sc Lemma 5.2.} {\it The reproducing kernel of the space
$W^a_{b,c}$ is  }
$$
K^a_{b,c}(z,u)=\frac{\Gamma(a+b)\Gamma(a+c)}{\Gamma(b+c)}
   \F\left[\begin{array}{c} a+b,\, a+c\\ b+c\end{array}; z\overline u\right]
.$$

{\sc  Proof.} Obvious.
Indeed, assume $\phi_w(z)=K^a_{b,c}(z,w)$.
Then for any holomorphic function
$f(z)=\sum_{k=0}^\infty c_kz^k$ we have
\begin{align*}
\langle f,\phi_w\rangle_{W^a_{b,c}}=  \sum_{k=0}^\infty c_k w^k
 \frac {\Gamma(a+b+k)\Gamma(a+c+k)}
        {k!\, \Gamma(b+c+k)}  \langle z^k, \, z^k\rangle =
\\=\sum_{k=0}^\infty c_k w^k= f(w)
.\end{align*}

{\sc Remark.} Kernels $K^a_{b,c}(z,u)$
are special cases of the kernels introduced by Gross, Richards
\cite{GrR} and also special cases of the kernels considered by
Odzijewicz
\cite{Odz}.

\smallskip

{\bf 5.3. Plancherel formula.}
Let $g\in W^a_{b,c}$. Consider the transform
\begin{multline*}
J^a_{b,c}g(s)=\frac{1}{|\Gamma(a+is)|^2\Gamma(b+c)}
  \inttt (1+x)^{-a-b}g\left(\frac x{x+1}\right)
\times \\ \times
             \FBC\densityx dx
.\end{multline*}

{\sc Theorem 5.3} {\it The operator $J^a_{b,c}$
is a unitary operator}
$$W^a_{b,c} \to L^2\left( \R_+,\densityabc\right).$$

{\sc Proof.} The image of the function $f(z)=z^k$ under
the transform $J^a_{b,c}$ is
\begin{equation}
\frac{S_k(s^2; a,b,c)}
       {\Gamma(a+b+k)\Gamma(a+c+k)}
,\end{equation}
Thus the required statement follows from (5.2)
and the orthogonality relations for the Hahn polynomials.

\smallskip

{\bf 5.4. Operator $d/dz$.} Denote by $W^\infty$
the space of functions holomorphic in the disk
$|z|<1$ and smooth up to the boundary.
Let $a>1$. By Lemma 1.1, for $f\in W^\infty$ the function
$J^a_{b,c}f$ is holomorphic in the strip $|\Im s|<1+\epsilon$.

\smallskip

{\sc Proposition 5.4.} {\it For each $f\in W^\infty$ and each $a>1$
$$J^a_{b,c}z\frac d{dz}f ={\cal L} J^a_{b,c} f ,$$
where $\cal L$ is the difference operator defined by} (4.1).

\smallskip

{\sc Proof.}
The image of the function $z^k$ is given by (5.3),
and the Hahn polynomials satisfy the difference equations (4.2).

\bigskip

{\large\bf \S 6. $\Lambda$-function and its properties}

\addtocounter{sec}{1}
\setcounter{equation}{0}

\nopagebreak

\bigskip

Motivation of the definition of the $\Lambda$-function
introduced in this Section is contained in above in 3.7.

{\bf 6.1. Definition.} Let $a,b,c\in \R$.
Let $x\in\R_+$.
We define the $\Lambda$-function by
\begin{multline}
\Lambda^a_{b,c}(x)=\frac 1{\pi\Gamma(b+c)} \inttt
\Gamma(a+is)\frac{\Gamma(b+is) \Gamma(b-is) \Gamma(c+is) \Gamma(c-is)}
      {\Gamma(2is)\Gamma(-2is)} \times\\ \times
\FBC\,
ds
.\end{multline}

{\sc Remark.}
The function $\Lambda^a_{b,c}(x)$ admits
a holomorphic continuation to the domain
 $\Re a, \Re b, \Re c>0$. It will be more pleasant for us
to formulate properties of the $\Lambda$-function
for real $a,b,c$.

\smallskip

I couldn't express the integral
(6.1) in terms of standard special functions
 (see the list in \cite{PBM})
by a finite number of algebraic operations
(except some special values of $b,c$,
see below 6.6).
I think that this is impossible.
The integrand in (6.1) is slightly intricate,
but the true reason of nontriviality of the integral (6.1)
are unusual limits of the integration.
Indeed,   let us  change
the limits of integration (for real $a,b,c$)  to
\begin{equation}
\int_{-\infty}^{+\infty}
.\end{equation}
This is equivalent to an evaluation of
\begin{equation}
\Re     \Lambda^a_{b,c}(x)
.\end{equation}
Then our integral becomes a Slater type integral
(see \cite{Sla}, \cite{Mar}) and the  standard residue  machinery
allows to obtain some expansions of
$\Re     \Lambda^a_{b,c}(x)$
into series.
Nevertheless I couldn't find a nice final expression for
 $\Re\Lambda^a_{b,c}(x)$
using this approach.

For some special values of $b,c$
 expression (6.3) can be evaluated explicitly (see 6.7)

Substituting  $x=0$ to (6.3), we obtain the Barnes integral
$$\Re     \Lambda^a_{b,c}(0)=
\frac 1{2\pi\Gamma(b+c)}\int_{-\infty}^{+\infty}
\frac{\Gamma(a+is)\Gamma(b+is) \Gamma(b-is) \Gamma(c+is) \Gamma(c-is)}
      {\Gamma(2is)\Gamma(-2is)} ds.$$
Applying the standard Barnes method,
 (see \cite{Sla}, \cite{Mar}, \cite {AAR})
we represent this integral as a sum of two hypergeometric series
${}_4 F_4(1)$
\begin{multline*}
\Re     \Lambda^a_{b,c}(0)=
\sum_0^\infty
\frac{(-1)^n \Gamma(a-b-n)\Gamma(2b+n)\Gamma(c-b-n)\Gamma(c+b+n)}
         {n!\,\Gamma(2b+2n)\Gamma(-2b-2n)}+\\+
\sum_0^\infty
\frac{(-1)^n  \Gamma(a-c-n)\Gamma(2c+n)\Gamma(b-c-n)\Gamma(b+c+n)}
           {n!\,\Gamma(2c+2n)\Gamma(-2c-2n)}
.\end{multline*}
We observe that $\Lambda^a_{b,c}(0)$ is some kind of a Barnes
integral over a nonclosed contour. Nevertheless
the expression $\Lambda^a_{b,c}(0)$
is 'better' than the 'indefinite Barnes integral'
\begin{equation}
\frac 1{2\pi\Gamma(b+c)}\int_{u}^{+\infty}
\Gamma(a+is)\frac{\Gamma(b+is) \Gamma(b-is) \Gamma(c+is) \Gamma(c-is)}
      {\Gamma(2is)\Gamma(-2is)} ds
,\end{equation}
since the point $s=0$ is a distinguished point for
the second factor of the integrand.

\smallskip

{\bf 6.2. Imitation of \O rsted problem.}
Consider the space $W^a_{b,c}$. The operator
$J^a_{b,c}$ is a unitary operator from $W^a_{b,c}$ to
the space
\begin{equation}
L^2\left(\R_+, \densityabc ds\right)
.\end{equation}
The multiplication operator
\begin{equation}
Mg(s)=\Gamma(a+is)g(s)
\end{equation}
is a unitary operator from the space (6.5) to the space
$$
L^2\left(\R_+, \densitybc ds\right)
$$
and the inverse index transform
$(J_{b,c})^{-1}$ is a unitary operator from (6.6) to
$$L^2(\R_+, x^{b+c-1}(1+x)^{b-c}dx).$$
Thus we obtain the unitary map
$$(J_{b,c})^{-1}\circ M \circ J_{b,c}^{a}: \,\,
  W_{b,c}^a\to   L^2(\R_+, x^{b+c-1}(1+x)^{b-c}dx)  .$$
This map is an imitation of the operator
$$V_\theta\to L^2(\B_m(\K)) $$
discussed above in 3.7.
The function $\Lambda^a_{b,c}$ is the image of the function $f(z)=1$
under this map.

\smallskip

{\bf 6.3. Direct corollaries of definition of
the $\Lambda$-function.}
The inversion formula for the index transform implies
$$
\frac 1{\Gamma(b+c)}\intt
\Lambda^a_{b,c}(x)\FBC\densityx dx=\Gamma(a+is)
.$$
The Plancherel formula for the inverse
index transform implies
\begin{multline*}
\intt
\Lambda^a_{b,c}(x)   \overline{\Lambda^a_{b,c}(x)}\densityx dx=\\=
\frac 1\pi\intt\densityabc ds
=\Gamma(a+b)\Gamma(a+c)\Gamma(b+c)
.\end{multline*}
(for a deduction of the last row see \cite{AAR}, 3.6).

\smallskip

{\bf 6.4. Differential--difference equations.}
The identity
$$\F(\alpha,\beta;\gamma;-x)=(1+x)^{\gamma-\alpha-\beta}
                    \F(\gamma-\alpha,\gamma-\beta;\gamma;-x)$$
(see \cite{HTF}, vol. 1, 2.1(23)) implies
 \begin{equation}
\Lambda^a_{b,c}(x)=(1+x)^{c-b}\Lambda^a_{c,b}(x)
.\end{equation}
Differentiating the integral (6.1) in the parameter $x$
by the formula
$$\frac d {dx}\, \F(\alpha,\beta;\gamma;x) =
  \frac{\alpha\beta}{\gamma}\F(\alpha+1,\beta+1;\gamma+1;x)
,$$
we obtain
\begin{multline*}
\frac d {dx} \Lambda^a_{b,c}(x)=
-\frac 1 {\Gamma(b+c+1)}\inttt \Gamma(a+is)
 \left|\frac{\Gamma(b+1+is)\Gamma(c+is)}{\Gamma(2is)}\right|^2
\times\\ \times
\F(b+1+is,b+1-is;b+c+1;-x)\,ds
=- \Lambda^a_{b+1,c}(x)
\end{multline*}
and finally
\begin{equation}
\frac d {dx} \Lambda^a_{b,c}(x)=-\Lambda^a_{b+1,c}(x)
.\end{equation}
Observe the following corollary from formulas (6.7), (6.8)
\begin{equation}
\Lambda^a_{b+k,c+l}(x)=(-1)^{k+l} (1+x)^{c+l-b-k}
\frac {d^l} {dx^l}(1+x)^{b+k-c}\frac {d^k} {dx^k}\Lambda^a_{b,c}(x)
.\end{equation}

By the formula (see \cite{AAR}, 2.5.7)
$$\Bigl( x \frac d {dx} +\gamma-1\Bigr) \F(\alpha, \beta;\gamma;x)=
              (\gamma-1)\F(\alpha, \beta;\gamma-1;x)
,$$
we obtain
\begin{multline*}
\bigl( x \frac d {dx} +b+c-1\bigr)
\Lambda^a_{b,c}(x)
=\\=
-\frac 1 {\Gamma(b+c-1)}\inttt \Gamma(a+is)
 \left|\frac{\Gamma(b+is)\Gamma(c-1+is)}{\Gamma(2is)}\right|^2
(c+is)(c-is)
\times\\ \times
\F(b+is,b-is;b+(c-1);-x)\,ds
.\end{multline*}
Representing $(c+is)(c-is)$ in the form
$$
(c+is)(c-is)=(c^2-a^2)+(2a+1)(a+is)-(a+is) (a+1+is)
,$$
we obtain
\begin{multline}
\bigl( x \frac d {dx} +b+c-1\bigr)
\Lambda^a_{b,c}(x)
=\\=
(c^2-a^2)\Lambda^a_{b,c}(x)  +
 (2a+1) \Lambda^{a+1}_{b,c}(x) - \Lambda^{a+2}_{b,c}(x)
\end{multline}
or
\begin{multline}
-x \Lambda^a_{b+1,c}(x)=
(c^2-a^2-b-c+1)\Lambda^a_{b,c}(x)+
 (2a+1) \Lambda^{a+1}_{b,c}(x) - \Lambda^{a+2}_{b,c}(x)
.\end{multline}

Certainly our list of difference equations is not complete.

\smallskip

{\bf 6.5. One integral.}

{\sc Proposition 6.1.} {\it For  $n=0,1,2,\dots$ }
\begin{align}
&2\Re \intt \Lambda^{a}_{b,c}(x)
\overline{\Lambda^{a+n}_{b,c}(x)}\densityx dx
=\nonumber \\
&=
(2a)_{2n} \Gamma(a+c)\Gamma(b+c)\Gamma(a+b)
\,\,{}_4 F_3\left[
\begin{array}{c} -n/2,-(n+1)/2, a+b, a+c\\
    -n+1, 2a, 2a+1 \end{array}; 1 \right]
.\end{align}

The expression ${}_4F_3$ in our case is not a series but a finite sum.
The right part of (6.12) also can be represented in the form
\begin{multline*}
\Gamma(a+b)\Gamma(a+c)\Gamma(b+c)(2a)_{2n}-  \\  -
n\Gamma(a+b+1)\Gamma(a+c+1)\Gamma(b+c)(2a+2)_{2n-2}+\\
 +n\Gamma(b+c)\sum_{k=2}^{[(n+1)/2]}
\frac{(n-k-1)\dots(n-2k+1)}{k!}
\Gamma(a+b+k)\Gamma(a+c+k)(2a+2k)_{2n-2k}
.\end{multline*}

{\sc Proof.}
By the Plancherel formula, our integral equals
\begin{align}
&\frac 2 \pi\intt \Re\bigl[\Gamma(a+is)\overline{\Gamma(a+is+n)}\bigr]
  \densitybc ds=\nonumber\\=
&\frac 1 \pi\intt \bigl[ (a+is)_n+ (a-is)_n \bigr] \densityabc ds
.\end{align}

{\sc Lemma 6.2.}
\begin{align*}
&(a+is)_n+ (a-is)_n=   \\=
&(2a)_{2n}- n(2a+2)_{2n-2} (a+is)(a-is)+\\
&\qquad + n\sum_{k=2}^{[(n+1)/2]}
\frac{ (-1)^k(-n+k+1)_{k-1}}{k!}(2a+2k)_{2n-2k}(a+is)_k(a-is)_k=\\=
& (2a)_{2n} \sum_{j=0}^{[(n+1)/2]} \frac{
(-n/2)_j(-(n+1)/2)_j(a-is)_j(a+is)_j}
  {j!\, (-n+1)_j (2a)_{2j} }
.\end{align*}

By  Lemma 6.2, we write  integral (6.13) in the form
\begin{multline*}
(2a)_{2n} \sum_{k=0}^{[(n+1)/2]}   \frac{
(-n/2)_j(-(n+1)/2)_j}
  {j!\, (-n+1)_j (2a)_{2j} }   \times\\ \times
\intt\left|\frac{\Gamma(a+j+is)\Gamma(b+is)\Gamma(c+is) }
              {\Gamma(2is)} \right|^2 ds=
\end{multline*}
$$=
(2a)_{2n} \sum_{k=0}^{[(n+1)/2]}   \frac{
(-n/2)_j(-(n+1)/2)_j}
  {j!\, (-n+1)_j (2a)_{2j} }
\Gamma(b+c)\Gamma(a+j+b)  \Gamma(a+j+c)
.$$
This implies the required result.

{\bf 6.6. Relations with the  $\lambda$-function.}
The function $\lambda(z,a)$ is given by
$$\lambda(z,a):=\int_0^a z^{-t}\Gamma(t+1)\,dt.$$
It will be more convenient for us to change notations
and to consider the function
$$
\lambda^*(z,b)=\int_b^{i\infty}z^{-t}\Gamma(t+1)\,dt$$
(it is the same indefinite integral with
another point of the origin).
The function $\lambda$
is an element of a relatively exotic family of special functions
(sometimes they are called the {\it Volterra type functions}).
This family includes also the functions
\begin{align*}
&\mu(z,p)=\intt \frac{t^p z^t dt}{\Gamma(t+1)};
&\mu(z,p,\rho)=\intt \frac{t^p z^{t+\rho} dt}{\Gamma(\rho+t+1)}; \\
&\nu(z)=\intt \frac {z^t dt}{\Gamma(t+1)};
&\nu(z, \rho)=\intt \frac {z^{t+\rho} dt}{\Gamma(\rho+t+1)};
.\end{align*}
Observe some superficial resemblance
of these functions with (6.1). Indeed, these integrals
are similar to the Barnes integrals, but we have 'incorrect'
limits of the integration.
The functions
$\mu(z,p)$, $\mu(z,p,\rho)$ were appeared in the Volterra work
(\cite{Vol}) on the fractional derivatives with logarithmic terms
(see also \cite{SKM}). The functions $\mu$, $\nu$
were intensively discussed in French and Belgian mathematical
journals in the first half of 1940ies
(Humbert, Poli, Colombo, Parodi and others, for instance, see,
 \cite{HP}, \cite{Col}).
The first mention on the $\lambda$-function
that I find is contained in the tables of McLachlan, Humbert, Poli
 \cite{MHP}.

A theory of the functions $\nu$, $\mu$
is contained in 'Higher transcendent functions'
\cite{HTF},
18.3 (the book contains also a large bibliography on this subject),
The integrals with the functions $\nu$, $\mu$, $\lambda$
 can be find in the tables of Prudnikov, Brychkov, Marichev,
 see corresponding subsections in volumes 3--5.

 Let us show that the functions
$$
\Lambda^a_{1/2+k,l}(x),\qquad \Lambda^a_{k,1/2+l}(x)
$$
can be expressed in terms of $\lambda^*(z,b)$
(in particular, this is valid for the $\Lambda$-functions
related with groups of the subseries $\OO(1,2n+1)$).
For instance, consider  the function
\begin{align*}
&\Lambda^a_{1/2,1}(x):=\\
&= \frac 1{\pi\Gamma(3/2)}\inttt
                \Gamma(a+is)\left|\frac{\Gamma(1/2+is)\Gamma(1+is)}
                     {\Gamma(2is)}\right|^2
        \frac 1 {2s \sqrt x} \sin(2s \arcsh \sqrt x)
           ds=  \\
&=-\frac {4}{\sqrt {\pi x}}\inttt \Gamma(a+is)is
 [e^{2is \arcsh \sqrt x} - e^{-2is \arcsh \sqrt x}]\,ds=\\
&=-\frac {4}{\sqrt {\pi x}}\inttt
(\Gamma(a+1+is) -a\Gamma(a+is))\times\\
&\times
 [(\sqrt{x^2+1}+x)^{2is}-(\sqrt{x^2+1}-x)^{2is}  ]\,ds
.\end{align*}
We obtained a sum of 4 integrals that can be expressed
by the function
$\lambda^*$.

Using formula (6.9), we obtain an expression of
$\Lambda^a_{1/2+k,1+l}(x)$ in terms of
the derivatives of $\lambda^*(x)$.
The case $\Lambda^a_{k,1/2+l}$ is similar.

 Observe that the derivatives of the function $\lambda^*$
can be expressed algebraically by the same function
 $\lambda^*(x)$
\begin{multline*}
\frac d {dz}\lambda^*(z,b)=-\int_b^{+i\infty} \frac t z z^{-t}
             \Gamma(t+1)dt=\\
= -\frac 1 z\int_b^{+i\infty} (\Gamma(t+2)-\Gamma(t+1))z^{-t} dt=
- \lambda^*(z,b+1)+\frac 1 {z} \lambda^*(z,b)
.\end{multline*}

\smallskip

{\bf 6.7. The cases of explicit evaluation for $\Re\Lambda^a_{b,c}$.}
Let $k,l=0,1,2,\dots$.
Let us show that the following functions
can be evaluated explicitly
$$\Re\Lambda^a_{1/2+k,l}(x),\qquad \Re\Lambda^a_{k,1/2+l}(x).$$
Let us calculate
\begin{align*}
&\Re \Lambda^{a}_{1/2,1}:=\const\int_{-\infty}^{+\infty}
\Gamma(a+is)\, s^2 \frac 1 {2s\sqrt x} \sin(2s \arcsh \sqrt x) ds=\\
&=\const
\int_{-\infty}^{+\infty}
(\Gamma(a+1+is) -a\Gamma(a+is))\times\\
&\times
 [(\sqrt{x^2+1}+x)^{2is}-(\sqrt{x^2+1}-x)^{2is}  ]\,ds
.\end{align*}
Using the formula (see, for instance, \cite{AAR}, (2.4.1))
$$\frac 1 {2\pi i} \int_{-i\infty}^{+i\infty}
  \Gamma(u+t) z^t dt = z^u e^{-z},$$
we obtain
\begin{multline*}
\const\Bigl[
(x+\sqrt{x^2+1})^{2a}((x+\sqrt{x^2+1})^2-a)
  \exp\bigl\{- (x+\sqrt{x^2+1})^2\bigr\}-
\\-
(\sqrt{x^2+1}-x)^{2a}((\sqrt{x^2+1}-x)^2-a)
   \exp\bigl\{- (\sqrt{x^2+1}-x)^2\bigr\}
\Bigr]
.\end{multline*}

To obtain an explicit expression for $\Re\Lambda^a_{1/2+k,l}(x)$,
we can apply formula (6.9).

{\sc Remark.} Slight inconvenience of the formulas of this
and the previous subsection disappears after the substitution
$$x=\sh^2t.$$

{\bf 6.8. Comments.} At the end of \S 3, I emphasized that there is
some arbitrariness in the choice of the
$\Lambda$-function. Some of properties of our list survive
for arbitrary choice of the function
$\psi$ (см. (3.10)), i.e. for  general integrals
having the form
\begin{multline*}
\frac 1{\pi\Gamma(b+c)} \inttt
\Gamma(a+is)e^{iu(s)}\frac{\Gamma(b+is) \Gamma(b-is) \Gamma(c+is) \Gamma(c-is)}
      {\Gamma(2is)\Gamma(-2is)} \times\\ \times
\FBC\,
ds
,\end{multline*}
where $\Im u(s)=0$.
This is valid for integrals from Subsection 6.3,
for differential-difference equations (6.7)-(6.9),
and for Theorem 7.1, Proposition 7.2, Theorem 8.1.

It seems important, that the $\Lambda$-function has
a collection of properties
that are not corollaries of the motivation given in
3.7, these properties are valid only for
 $\psi=1$ (or $u=0$).
This is the difference equation (6.10)
and all material of 6.5-6.7.

It should be noted that the differential-difference equations
obtained in 6.4 and the expressions
by the $\lambda$-function are valid also for the indefinite
Slater integral (6.4).
All other properties of the $\Lambda$-function obtained above
don't survive  in this generality.

\bigskip

{\large\bf \S 7. Bireflected bases and generalized translate hypergroup}

\nopagebreak

\bigskip

\addtocounter{sec}{1}
\setcounter{equation}{0}

{\bf 7.1. Bireflected bases.}
In 5.3 we defined the unitary operator
$$Q=(J_{b,c})^{-1}\circ M \circ J_{b,c}^{a}: \,\,
  W_{b,c}^a\to   L^2(\R_+, x^{b+c-1}(1+x)^{b-c}dx)  .$$
By
$$\Xi_n(x)=\Xi_n(x;a;b,c)$$
we denote the image of the function $z^n$ under the transform
$Q$.

\smallskip

{\sc Theorem 7.1.}   a) {\it
The system of the functions
$\Xi_n$ forms an orthogonal basis in  $L^2(\R_+, x^{b+c-1}(1+x)^{b-c}dx)$.}

\smallskip

b)  {\it The functions $\Xi_n$ can be expressed in terms
of the $\Lambda$-function
by}
\begin{multline*}
\Xi_n(x;a,b,c)=\\
=\frac{(b+c)_n}{\Gamma(a+c)\Gamma(a+b+n)}
 \sum_{j=0}^n \frac{(-1)^j(-n)_j}{j!\, (a+c)_j (b+c)_j}
   x^{-b-c+1} \frac {d^j}{dx^j} x^{b+c+j-1} \Lambda^a_{b,c+j}
.\end{multline*}

\smallskip

{\sc Proof.}  The statement a) is obvious and we must
prove only b). Recall that the image of the function $z^n$
under the operator
$J_{b,c}^{a}$ is a Hahn polynomial.
Thus it is sufficient to evaluate the integral
\begin{multline*}
\frac 1{\Gamma(a+b+n) \Gamma(a+c+n)     \Gamma(b+c)}
\inttt\Gamma(a+is) \densitybc \times
\\ \times \FBC S_n(x^2; a,b,c)\,ds
.\end{multline*}
By the symmetry of the Hahn polynomials in $a,b,c$,
we represent the integral in the form
\begin{multline*}
\frac 1{\Gamma(a+b+n) \Gamma(a+c+n)     \Gamma(b+c)}
\inttt\Gamma(a+is) \densitybc \times
\\ \times \FBC  (a+c)_n(b+c)_n \FF
\left( \begin{array}{c} -n,\, c+is,\, c-is\\ a+c,\, b+c\end{array};1\right)
\,ds
.\end{multline*}

Expanding $\FF$ into the sum, we obtain
\begin{multline*}
\frac{(b+c)_n} {\Gamma(a+b+n) \Gamma(a+c)\Gamma(b+c)}
  \sum_{j=0}^n \frac{(-n)_j} {j! (a+c)_j (b+c)_j} \times
\\ \times
\inttt  \Gamma(a+is)
 \bigl|\frac{\Gamma(b+is)\Gamma(c+j+is)}{\Gamma(2is)}\bigr|^2 \FBC\,ds
.\end{multline*}
Applying the formula (\cite{HTF}, vol.1, 2.8(22))
$$(\mu)_k y^{\mu-1} \F(\alpha,\beta;\mu;y)=
\frac {d^k}{dy^k} y^{\mu+k-1}\F(\alpha,\beta;\mu+k;y)
,$$
we obtain
\begin{multline*}
\frac{(b+c)_n} {\Gamma(a+b) \Gamma(a+c)}
\sum_{j=0}^n \frac{(-n)_j} {j! (a+c)_j (b+c)_j}
\frac 1 {\Gamma(b+c+j)}
(-1)^j x^{1-b-c}\frac {d^j} {dx^j} x^{b+c+j-1}
\times \\ \times
 \inttt
  \Gamma(a+is) \left|\frac{\Gamma(b+is)\Gamma(c+j+is)}{\Gamma(2is)}\right|^2
          \F(b+is, b-is; b+(c+j);-x)\, ds
.\end{multline*}
This gives the required expression.

\smallskip

{\bf 7.2. Generalized translate operator.}
The generalized translate hypergroup
is an essential element of the theory of the index hypergeometric transform.
(see \cite{FK}, \cite{Koo2}). We only will slightly touch
on this subject.

Consider the set $J\subset \R_+\times \R_+ \times \R_+$
consisting of all points $(x,y,z)$ such that three numbers
$\arcsh \sqrt x, \arcsh \sqrt y, \arcsh \sqrt z$
satisfy the triangle inequality.
Consider the function $K(x,y,z)$ on $\R_+\times \R_+ \times \R_+$
that is identical zero outside the set $J$ and is equal to
\begin{multline*}
K(x,y,z)=
 \frac{2^{-4b}\Gamma(b+c)}{\Gamma(b+c-1/2)}
\frac{ \bigl( (1+x)(1+y)(1+z)\bigr)^{c-1}}
    {(xyz)^{b+c-1}} \times
\\ \times(1-B^2)^{b+c-1/2}
 \F(2b-1, 2c-1; b+c-1/2; (1-B)/2)
\end{multline*}
on the set $J$, where $B$ denotes the expression
$$B=\frac{x+y+z+2}{2(1+x)(1+y)(1+z)}.$$
The {\it generalized translate operator}
$T_z$ в $L^2(\R_+,\densityx)$ is given by
$$T_z f(y)=\intt K(x,y,z) f(x)\densityx dx.$$

It is known that
\begin{equation}
 J_{b,c} T_z (J_{b,c})^{-1} g(s)= \F(b+is,b-is;b+c;-z) f(s)
.\end{equation}

{\sc Proposition 7.1.}
\begin{multline*}
\inttt \inttt \Lambda^a_{b,c}(x)
\overline{\Lambda^a_{b,c}(x) }
K(x,y,z)\,\,\densityx y^{b+c-1}(1+y)^{b-c} dx\,dy
=\\
=\Gamma(a+b)\Gamma(a+c)\Gamma(b+c)(1+z)^{-a-b}
.\end{multline*}

{\sc Proof.} The expression in the left side is the scalar
product
$\Lambda^a_{b,c}$ с $T_z\Lambda^a_{b,c}$ in $L^2(\R_+,\densityx)$.
By the Plancherel formula (7.1), this equals
\begin{multline*}
 \intt \densitybc
 \Gamma(a+is)\overline{\Gamma(a+is)}\F(b+is,b-is; b+c;-z)\,ds
.\end{multline*}
The last expression is the inverse index transform
of $|\Gamma(a+is)|^2$. But it was evaluated in Lemma 4.1.

\smallskip

{\bf 7.3. Generalized translate operator and bireflected bases.}
It seems that the generalized translate operator has a simple explicit
matrix in the bireflected basis.

Author can evaluate only the first row.

\smallskip

{\sc Proposition 7.3.}
$$T_z  \Lambda^a_{b,c} = T_z \Xi_0(x; a,b,c)=
  (1+z)^{a-b}\sum_{n=0}^\infty \frac 1 {n! (b+c)_n}
\Bigl( \frac z {z+1}\Bigr)^n
          \Xi_n(x; a,b,c).$$

{\sc Proof.}
The index transform of
$T_z  \Lambda^a_{b,c}$ is
\begin{align*}
&\Gamma(a+is)\F(b+is,b-is; b+c;-z)=
 \\
&=\Gamma(a+is)(1+z)^{-b-is}\F(b+is,c+is;b+c; \frac z {z+1})=
\\
&= \Gamma(a+is) (1+z)^{a-b}\Bigl(1-\frac z {z+1}\Bigr)^{a+is}
\F(b+is,c+is;b+c; \frac z {z+1})
.\end{align*}
Using the following generating function
for the continuous dual Hahn polynomials
(see \cite{AAR}, p.349),
we obtain
$$
\Gamma(a+is) (1+z)^{a-b}
\sum_{k=0}^\infty
\frac{S_k(x^2; a,b,c)}{k! \, (b+c)_k} \left( \frac z {z+1}\right)^k.$$
This is equivalent to the required statement.

\bigskip

{\large\bf \S 8. An application of $\Lambda$-function and
of bireflected bases.
Symmetric spaces $\U(p,q)/\U(p)\times \U(q)$}

\nopagebreak

\bigskip

\addtocounter{sec}{1}
\setcounter{equation}{0}

{\bf 8.1. The space $\U(p,q)/\U(p)\times\U(q)$.}
Assume $p\le q$.
Denote by $\U(p,q)$ the pseudounitary group of the order $(p,q)$,
i.e. the group of all matrices
$g=
\bigl(\begin{smallmatrix}a&b\\c&d\end{smallmatrix}\bigr)$,
satisfying the condition
$$
\begin{pmatrix}a&b\\c&d\end{pmatrix}
\begin{pmatrix}1&0\\0&-1\end{pmatrix}
\begin{pmatrix}a&b\\c&d\end{pmatrix}^*=
\begin{pmatrix}1&0\\0&-1\end{pmatrix}
.$$

Denote by $\B_{p,q}$ the space of all matrices $z$ over $\C$
having the size $p\times q$ with the norm satisfying $\|z\|<1$.
The group $\U(p,q)$ acts on $\B_{p,q}$ by the linear fractional
transforms
$$
z\mapsto z^{[g]}:= (a+zc)^{-1}(b+zd)
.$$

The stabilizer $K$ of the point 0 consists of matrices having the form
$$
\begin{pmatrix}a&0\\0&d\end{pmatrix}
  ;\qquad a\in\U(p), d\in\U(q)
.$$
The space $\B_{p,q}$ is a homogeneous symmetric space,
$$\B_{p,q}=\U(p,q)/\U(p)\times\U(q).$$
The  $\U(p,q)$-invariant measure on $\B_{p,q}$ has the density
$\det(1-zz^*)^{-p-q}$ with respect to the Lebesgue measure.
The group $\U(p,q)$ acts in the space $L^2$
with respect to the invariant measure
by the substitutions
\begin{equation}
T(g)f(z)=f(z^{[g]})
.\end{equation}

\smallskip

{\bf 8.2. Berezin spaces $H^\circ_\theta(\B_{p,q})$.}
The Berezin space $V^\circ_\theta(\B_{p,q})$
is  the space of functions on
$\B_{p,q}$ defined by the kernel
$$K(z,u)=\frac{\det(1-zz^*)^{\theta/2} \det(1-uu^*)^{\theta/2}}
             {|\det(1-zu^*)|^\theta}
.$$
The group  $\U(p,q)$ acts in $V^\circ_\theta(\B_{p,q})$  by formula
(8.1).
For sufficiently large $\theta$ the representation of
$\U(p,q)$ in $V^\circ_\theta(\B_{p,q})$
is equivalent to the representation in  $L^2(\B_{p,q})$
(see \cite{Ber2}, see also \cite{OO}, \cite{Ner3}).
The {\it \O rsted problem} is to construct
explicitly a unitary intertwining operator between
these representations.

\smallskip

{\bf 8.3. $\Lambda$-function of the space $\U(p,q)/\U(p)\times\U(q)$.}
Denote by 
$\lambda_1(z),\dots,\lambda_p(z)$
the singular values of the matrix $z$. Define
the new variables by
$$x_j(z)=\lambda_j^2(z)/(1-\lambda_j^2(z)).$$
We define the $\Lambda$-function of the space $\U(p,q)/\U(p)\times\U(q)$
by
\begin{equation}
\Lambda^\theta(z)=\det\begin{pmatrix}
  \Xi_0(x_1)&\dots
          & \Xi_0(x_p)\\
  \Xi_1(x_1)&\dots
          & \Xi_1(x_p)\\
\vdots & \ddots & \vdots \\
  \Xi_{p-1}(x_1)&\dots
          & \Xi_{p-1}(x_p)
\end{pmatrix}
,\end{equation}
where
$$\Xi_j(x)=\Xi_j(x; \theta-(q+p-1)/2; (q-p+1)/2;(q-p+1)/2) $$
are the first $p$ elements of the bireflected basis.

\smallskip

{\sc Theorem 8.1.} {\it For $g\in\U(p,q)$ by $u$ we denote
the image of the point $0\in\B_{p,q}$ under $g$. Then}
$$\int\limits_{\B_{p,q}}
 \Lambda^\theta(z)\overline{\Lambda^\theta(z^{[g]})}
  \det(1-zz^*)^{-p-q} dz=\det(1-uu^*)^{\theta/2}.$$

\smallskip

{\sc Corollary 8.2.} {\it Define the function $L(z,u)$
on $\B_{p,q}\times\B_{p,q}$ by the formula
$$L(z,u)= \Lambda(z^{[g_u]}),$$
where $g_u\in\U(p,q)$ is any element that transforms
the point
$0\in\B_{p,q}$ to $u$. Then the operator
$$Mf(z)=\int\limits_{\B_{p,q}}   L(z,u)f(u)\det(1-uu^*)^{-p-q}\,du$$
is a unitary $\U(p,q)$-intertwining operator}
$$L^2(\B_{p,q})\to V^\circ_\theta(\B_{p,q}).$$

\smallskip

{\sc Remark.} The definition of the kernel $L(z,u)$ can be reformulated
in the following form.
It is defined by the determinant
(8.2), where the numbers
$\lambda_j$ are the singular values of the matrix
$$(1-zz^*)^{-1/2}(1-zu^*)(1-uu^*)^{-1/2}.$$

\smallskip

{\sc Remark.} The similar statement is valid
for the groups of series
$\OO(n,\C)$ и $\Sp(2n,\C)$.

{\sc Institute of Theoretical  and Experimental Physics}

neretin@main.mccme.rssi.ru

\end{document}